\documentstyle[12pt]{article}

\begin{document}

\begin{titlepage}
\begin{center}
{\Large\bf Classical and quantum duality in jordanian quantizations 
}
\end{center}
\vskip 1cm
\centerline{\large \bf
P. P. Kulish\footnote{Supported by the RFFI grants NN 96-01-00851
and 98-01-00310.}}
\vskip 1mm
\begin{center}
{St.Petersburg
Department of the Steklov
Mathematical Institute,}\\
{Fontanka 27, St.Petersburg, 191011,
Russia }
\end{center}
\vskip 5mm
\centerline{\large \bf
V. D. Lyakhovsky\footnote{Supported by the RFFI grant N 97-01-01152 and by
the Direction General de Ensenanza Superior de 
la Ministerio de Educacion y Cultura of Spain, grant N SAB1995-0610}}
\vskip 1mm
\begin{center}
{Department of Physics}\\
{St.Petersburg
State University}\\
{Ulianovskaya 1, Petrodvorets, St.Petersburg, 198904,
Russia }
\end{center}
\vskip 5mm

\begin{abstract}
The limiting transitions between different types of quantizations are studied
by the deformation theory methods.
We prove that for the first order coboundary deformation of a Lie
bialgebra $\left( g,g_1^{*}\right) \longrightarrow \left( g,g_1^{*} + \xi
g_2^{*}\right)$ one can always get the quantized Lie bialgebra
${\cal A}\left( g,g_2^{*}\right)$ as a limit of the sequence of
quantizations of the type ${\cal A}\left( g,g_1^{*}\right)$. 
The obtained results are illustrated by some low-dimensional examples of 
quantized Lie algebras and superalgebras. 
\end{abstract}
\end{titlepage}

\section{Introduction}

The Drinfeld and Jimbo deformations of universal enveloping of
simple Lie algebras \cite{DRIN,JIMB} correspond to Lie bialgebras with
classical $r$-matrix
$$
r_{DJ}=\sum_{i=1}^kt_{ij}H_i\otimes H_j+\sum_{\alpha \in \Delta
_{+}}E_\alpha \otimes E_{-\alpha },
$$
where $k$ is the rank, $t_{ij}$ is the inverse Cartan matrix, and $\Delta
_{+}$ is the set of positive roots.
The triangular quantum groups and twistings \cite{D2}, for
example, the jordanian quantization of $sl(2)$ or of its Borel subalgebra $
{\bf B}_{+}$ ($\{h,x|[h,x]=2x\}$) with $r=h\otimes x-x\otimes h=h\wedge x$
\cite{DRIN}, have the triangular $R$-matrix ${\cal R}={\cal F}_{21}{\cal F}
^{-1}$ defined by the twisting element \cite{GER,OGIEV}
\begin{equation}
\label{og-twist}{\cal F}=\exp \{\frac 12h\otimes \ln (1+2\xi x)\}.
\end{equation}

The extension of this twist \cite{KLM} implies the existence of a
subalgebra ${\bf L}$ with special properties. This is a solvable subalgebra
with at least four generators. All simple Lie algebras except $sl(2)$
contain such ${\bf L}$ and in any of them a deformation induced by twist of $
{\bf L}$ can be performed. In particular for $
{\cal U}(sl(N))$ the following twisting element ${\cal F}\in {\cal U}
(sl(N))^{\otimes 2}$ can be applied,
\begin{equation}
\label{twist-sl(N)}{\cal F}=\exp \{2\xi \sum_{i=2}^{N-1}E_{1i}\otimes
E_{iN}e^{-\sigma }\}\exp \{H\otimes \sigma \},
\end{equation}
where $x=E_{1N}$, $H=E_{11}-E_{NN}$, $\sigma =\frac 12\ln (1+2\xi x)$.

In this work the cohomological interpretation of the interrelations
between the Drinfeld-Jimbo (or standard) quantum algebra ${\cal U}_q(sl(N))$
and the jordanian (or non-standard) ${\cal U}_\xi (sl(N))$ are discussed.
The existence of such interrelations was already pointed out in \cite{GER}.
The operator $\exp (\xi {\rm ad}E_{1N})$ (with the highest root generator 
$E_{1N})$
turns $r_{DJ}$ into the sum $r_{DJ}+\xi r_j$. Hence,
\begin{eqnarray}r_j=-\xi \left( H_{1N}\wedge
E_{1N}+2\sum_{k=2}^{N-1}E_{1k}\wedge
E_{kN}.\right) ,
\end{eqnarray}is a classical $r$-matrix too.

Here the decomposition of $r$-matrix can be treated as a consequence
of a specific similarity transformation. We shall show that the effect is
in fact connected with the quite general Lie bialgebra deformation
properties.

\section{Deformed coboundary Lie bialgebras}

Let us consider a coboundary Lie bialgebra $\left( g,g^{*}\right) $
and its first order deformation which we shall write in a most general form
: $\left( g_h,g_\xi ^{*}\right) $ . We can consider it as a four dimensional
variety
$$
\left( \alpha _1g_1+\alpha _2g_2,\beta _1g_1^{*}+\beta _2g_2^{*}\right)
=\left( g_{\alpha _1,\alpha _2},g_{\beta _1,\beta _2}^{*}\right) .
$$
Being the first order deformation this pair contains in fact four
Lie bialgebras:
$ \left( g_i,g_k^{*}\right) _{i,k=1,2},\quad {\rm with\
compositions\quad }\mu _i,\mu _k^{*}$.
As a coboundary Lie bialgebra the pair $\left( g_{\alpha _1,\alpha
_2},g_{\beta _1,\beta _2}^{*}\right) $ corresponds to the classical $r$
-matrix satisfying the MCYBE, its symmetric part being ad$^{\otimes 2}$
-invariant. Hence $g_{\beta _1,\beta _2}^{*}$ defines a Lie coalgebra
(on the space $V_g$ of the Lie algebra $g_{\alpha _1,\alpha _2}$ ),
\begin{equation}
\label{colie-gen}
\begin{array}{ccccc}
\delta _{\beta _1,\beta _2}\left( x\right) & = & \alpha _1\mu _1\left(
x\otimes 1+1\otimes x,r\right) & + & \alpha _2\mu _2\left( x\otimes
1+1\otimes x,r\right) \\
& = & \left( \beta _1\delta _1(x)\right) & + & \beta _2\delta _2\left(
x\right) .
\end{array}
\end{equation}
When the parameters $\alpha _i,\beta _k$ are independent
no solutions of (\ref{colie-gen}) smoothly depending on parameters can
be found.

We shall confine ourselves to the
following two types of restrictions on parameters
$\left\{ \alpha _i,\beta _k\right\} $
$$
I. \quad  \alpha _i=\beta _i, \quad \quad II. \quad  \alpha _2=0.
$$

In the first case the constant $r$-matrix can serve as a solution of (\ref
{colie-gen}). This situation is characteristic for the quantum double case
\cite{Lya}. The equation (\ref{colie-gen}) factorizes:%
$$
\delta _i\left( x\right) =\mu _i\left( x\otimes 1+1\otimes x,r\right) .
$$
Thus for a quantum double of two Borel algebras:%
$$
\begin{array}{ll}
 \left[ H,X_{\pm }\right] =\pm \alpha _1X_{\pm }, &
\Delta X_{+}=X_{+}\otimes 1+e^{-\theta \alpha _2H}\otimes X_{+}, \\
 \left[ H^{\prime },X_{\pm }\right] =\pm \alpha _2X_{\pm }, &
\Delta X_{-}=X_{-}\otimes e^{\theta \alpha _1H^{\prime }}+1\otimes X_{-},
\\
 \left[ X_{+},X_{-}\right] =\frac 1\theta \left( e^{\theta \alpha
_1H^{\prime }}-e^{-\theta \alpha _2H}\right); &
\end{array}
$$
the Lie bialgebra has the form $\left( \alpha _1g_1+\alpha _2g_2,\alpha
_1g_1^{*}+\alpha _2g_2^{*}\right) $ , where%
$$
\begin{array}{c}
g_1=\left\{
\begin{array}{c}
\left[ H,X_{\pm }\right] =\pm X_{\pm }; \\
\left[ X_{+},X_{-}\right] =H^{\prime };
\end{array}
\right\} ,\ g_2=\left\{
\begin{array}{c}
\left[ H^{\prime },X_{\pm }\right] =\pm X_{\pm }; \\
\left[ X_{+},X_{-}\right] =H;
\end{array}
\right\} , \\
g_1^{*}=\left\{ \left[ \widehat{H^{\prime }},
\widehat{X_{-}}\right] =-\theta \widehat{X_{-}};\right\} ,\
g_2^{*}=\left\{ \left[ \widehat{H},\widehat{X_{+}}\right] 
=-\theta \widehat{X_{+}};\right\} .
\end{array}
$$
and $\left\{\widehat{H},\widehat{H^{\prime }},\widehat{X_{\pm}}\right\}$ is
the dual base. 
One can easily see that in this situation
to be a solution of the equation (\ref{colie-gen}) the
classical $r$-matrix must be constant :%
$$
r=\theta \left( X_{+}\otimes X_{-}+H\otimes H^{\prime }\right).
$$

In the second case the equation can be solved only with $r$-matrix
depending on $\beta _i$ and, moreover, the latter can be decomposed
into the sum of two constituent $r$-matrices for $\left( g,\beta
_1g_{1}^{*}\right) $ and $\left( g,\beta g_{2}^{*}\right) $ correspondingly
$$
r=\beta _1 r_1+\beta _2 r_2.
$$

This latter case of the Lie bialgebras first order deformations
appear in the studies of jordanian quantizations. Below we shall treat them 
in details.

\section{Drinfeld-Jimbo and jordanian deformations as mutually first order
ones}

It is well known that some sorts of jordanian deformations can be
treated as limiting structures for certain sequences of standard
quantizations \cite{GER,OGIEV,ABD}. These properties are more transparent
for quantum groups.

Let the $N\times N$-matrix $T$ formed by the generators of the
standard (FRT-deformed) quantum group $Fun_h(SL(N))$ be subject
to the similarity transformation with the matrix
\begin{equation}
\label{transM}M=1+\frac \xi {q-1}\rho \left( E_{1N}\right) \quad (q=e^h)
\end{equation}
(the coproduct ($\Delta T=T\stackrel{.}{\otimes }T$) is conserved). As far
as $q\neq 1$ the transformed group $Fun_{h;\xi }(SL(N))$ is equivalent to
the original one. Compare the corresponding Lie bialgebras: $\left(
g,g_{h;0}^{*}\right) =\left( sl(N),\left( sl(N)\right) ^{*}\right) $ and $%
\left( g,g_{h;\xi }^{*}\right) $. Here $g=sl(N)$ is fixed and only the
second Lie multiplication law $\left( \mu _{h;0}^{*}:V_{g^{*}}\wedge
V_{g^{*}}\rightarrow V_{g^{*}}\right) $ changes:%
$$
\mu _{h;0}^{*}\rightarrow \mu _{h;\xi }^{*}.
$$
One can check that the new Lie product decomposes as:
\begin{equation}
\label{sum-mu'}\mu _{h;\xi }^{*}=\mu _{h;0}^{*}+\xi \mu ^{\prime }.
\end{equation}
$\mu ^{\prime }$ is defined by the transformed $RTT=TTR$ equations.
Change the
coordinate functions of $SL(N)$ arranged in $T$ for the exponential ones $
T=\exp (\epsilon Y)$ and the parameters $h\longmapsto \epsilon h,\ \xi
\longmapsto \epsilon \xi $ . Tending $\epsilon $ to zero one gets both
summands in (\ref{sum-mu'}). In the canonical $gl(N)$-basis the second
one is:
\begin{equation}
\label{dual-sln1}
\begin{array}{l}
\mu ^{\prime }\left( Y_{1k},Y_{ij}\right) =2\delta _{ik}Y_{Nj},
{\rm \ for}\ k,j<N;\quad i>1, \\ \mu ^{\prime }\left( Y_{ij},Y_{lN}\right)
=-2\delta _{jl}Y_{Nj},\
{\rm for}\ j<N;\quad i,l>1, \\ \mu ^{\prime }\left( Y_{ij},Y_{1N}\right)
=-\delta _{j1}Y_{i1}-\delta _{iN}Y_{Nj},\ {\rm for}\ j<N;\quad i>1,
\end{array}
\end{equation}
\begin{equation}
\label{dual-sln2}
\begin{array}{l}
\mu ^{\prime }\left( Y_{1i},Y_{1N}\right) =-Y_{1i},\
{\rm for}\ N>i>1, \\ \mu ^{\prime }\left( Y_{1N},Y_{kN}\right) =Y_{kN},\
{\rm for}\ k<N<1, \\ \mu ^{\prime }\left( Y_{11},Y_{1N}\right) =\mu ^{\prime
}\left( Y_{1N},Y_{NN}\right) =-\left( Y_{11}-Y_{NN}\right) , \\
\mu ^{\prime }\left( Y_{1i},Y_{1k}\right) =\delta _{i1}Y_{Nk},
{\rm \ for}\ k,i<N;\quad k>1, \\ \mu ^{\prime }\left( Y_{iN},Y_{kN}\right)
=-\delta _{kN}Y_{i1},
{\rm \ for}\ k,i>1;\quad i<N, \\ \mu ^{\prime }\left( Y_{1i},Y_{kN}\right)
=\delta _{i1}Y_{k1}-\delta _{kN}Y_{Ni}-2\delta _{ik}\left(
Y_{11}-Y_{NN}\right) ,\ {\rm for}\ i<N;\quad k>1.
\end{array}
\end{equation}
The map $\mu ^{\prime }$ not only defines the
infinitesimal deformation of $\mu _{h;0}^{*}$ but is itself a Lie
multiplication. The full map
$\mu _{h;\xi }^{*}$ is a deformation of $\mu _{h;0}^{*}$ and $%
\mu ^{\prime }$ does not depend on $h$ or $\xi :$
\begin{equation}
\label{sum-mu}
\mu _{h;\xi }^{*}=\mu _{h;0}^{*}+\mu _{0;\xi }^{*}.
\end{equation}
Thus the composition $\mu _{h;\xi }^{*}$ is a Lie multiplication deformed
in the first order.

Both summands are the Lie maps and at the same time can be considered as
deforming functions of each other:
$$
\mu _{h;0}^{*}\in Z^2\left( g_{0;\xi }^{*},g_{0;\xi }^{*}\right) , \quad
\mu _{0;\xi }^{*}\in Z^2\left( g_{h;0}^{*},g_{h;0}^{*}\right) .
$$
The equivalence
$Fun_{h;\xi }(SL(N))\approx $ $Fun_h(SL(N))$ (for $h\neq 0$) signifies that $
\mu _{0;\xi }^{*}$ is in fact a coboundary,%
$$
\mu _{0;\xi }^{*}\in B^2\left( g_{h;0}^{*},g_{h;0}^{*}\right) .
$$
On the contrary, the composition $\mu _{h;0}^{*}$ corresponds to a
nontrivial cohomology class%
$$
\mu _{h;0}^{*}\in H^2\left( g_{0;\xi }^{*},g_{0;\xi }^{*}\right) ,
$$
the deformation of $\mu _{0;\xi }^{*}$ by $\mu _{h;0}^{*}$ is essential \cite
{NIJ}.

With respect to the cochain complex $C^n:\bigwedge^nV_g\rightarrow
V_g\wedge V_g$, all the Lie algebras $g_{h;\xi }^{*}$ are dual to one and
the same $g=sl(N)$:
$$
\mu _{h;0}^{*},\mu _{0;\xi }^{*}\in B^1\left( g,g\wedge g\right) .
$$

Thus (see Section 1) the classical $r$-matrix of ${\cal U}%
_{h;\xi }\left( sl(N)\right) \approx $ $\left( Fun_{h;\xi }(SL(N))\right)
^{*}$ must exhibit the decomposition property:
\begin{equation}
\label{r-sln}
\begin{array}{c}
r_{h;\xi }=r_{h;0}+r_{0;\xi }=\frac hN\left( \sum_{k=1}^{N-1}k\left(
N-k\right) H_{k,k+1}\otimes H_{k,k+1}\right. \\
+\left. \sum_{k<l}\left( N-l\right) k\left( H_{k,k+1}\otimes
H_{l,l+1}+H_{l,l+1}\otimes H_{k,k+1}\right) \right) \\
+2h\sum_{k<l}\left( E_{lk}\otimes E_{kl}\right) \\
-\xi H_{1N}\wedge E_{1N}-2\xi \sum_{k=2}^{N-1}E_{1k}\wedge E_{kN}.
\end{array}
\end{equation}
In the limit $h\rightarrow 0$ one gets the element
\begin{eqnarray}
\lim _{h\rightarrow 0}r_{h;\xi }=r_{0;\xi }=-\xi
\left( H_{1N}\wedge
E_{1N}+2\sum_{k=2}^{N-1}E_{1k}\wedge E_{kN}.\right) ,
\label{CYBE0}
\end{eqnarray}
that coincides with the $r$-matrix that can be obtained from ${\cal R}$
\begin{equation}
\label{R-sl(3)}
\begin{array}{c}
{\cal R}={\cal F}_{21}{\cal F}^{-1} \\ =\prod_j\exp \left( 2\xi
E_{jN}e^{-\sigma }\otimes E_{1j}\right) \exp \left( \sigma \otimes
H_{1N}\right) \exp \left( -H_{1N}\otimes \sigma \right) \cdot \\
\cdot \prod_j\exp \left( -2\xi E_{1j}\otimes E_{jN}e^{-\sigma }\right)
\end{array}
\end{equation}
of the twisted algebra ${\cal U}_{{\cal F}}(sl(N))$ \cite{KLM}.

To clarify the contraction properties of $Fun_{h;\xi }(SL(N))$ let
us consider the 1-parameter subvariety $\left\{ g_{h;1-h}^{*}\right\} $ of
Lie algebras $g^{*}{}_{h;\xi }^{}$ (putting $\xi =1-h$ in (\ref{sum-mu})).
Each dual pair $\left( sl(N),g_{h;1-h}^{*}\right) $ is a Lie bialgebra and
thus is quantizable \cite{ETI}. The result is the set ${\cal A}_{s;h}$. This
set is smooth (in the formal series topology). The 1-dimensional
boundaries ${\cal A}_{0;h}$ and ${\cal A}_{s;0}$ of ${\cal A}_{s;h}$ are the
quantizations of $\left( sl(N),g_{1;0}^{*}\right) $ and of $\left(
sl(N),g_{0;1}^{*}\right) $ respectively. Each internal point
in ${\cal A}_{s;h}$ can be connected with a boundary by a smooth
parametric curve $a(u)$ starting in ${\cal A}_{0;h}$ and ending
in ${\cal A}_{s;0}$. This does not necessarily
mean that this limit is a faithful contraction -- it may not be in orbit.
This is just what happens in our case. For every positive $h$
fixed the subset $\left\{ Fun_{h;\xi }(SL(N))\right\} $ is in the $GL(N^2)$
-orbit of the corresponding $Fun_{h;0}(SL(N))$. To attain the
points $Fun_{0;\xi }(SL(N))$ one must tend $h$ to zero by crossing the set of
orbits.

One of the principle conclusions is that the possibility to obtain
the jordanian deformation $Fun_{0;\xi }(SL(N))$ as a limiting transformation
of the standard quantum group -- $Fun_{h;0}(SL(N))$ (and on the dual list to
get the twisted $q$-algebra ${\cal U}_{{\cal F}}(sl(N))$ as a limit of the
variety of standard quantized algebras ${\cal U}_q(sl(N)))$ is provided by
the fact that the 1-cocycle $\mu _{0;\xi }^{*}\in Z^1(sl(N),sl(N)\wedge
sl(N))$ (that characterizes the Lie bialgebra for ${\cal U}_{{\cal F}%
}(sl(N)) $ ) is at the same time the 2-coboundary $\mu _{0;\xi }^{*}\in
B^2\left( g_{h;0}^{*},g_{h;0}^{*}\right) $ the Lie algebra $g_{h;0}^{*}$
being the standard dual of $sl(N)$.

We want to stress that these interrelations of Lie maps are not specific
only to the factorizable twists such as (\ref{twist-sl(N)}).
Consider the $R$-matrix depending on two parameters ($h$ and $\xi$)
proposed in \cite{KUL} for the superalgebra ${\cal U}(osp(1|2))$.
It leads to a smooth sequence of deformations whose
limit (for $h \rightarrow 0$) was proved to be a twist.
The twisting element consists of two factors
$ {\cal F}^{(sj)} = {\cal F}^{(s)}{\cal F}$, where the first operator
${\cal F}$ is the jordanian twist (\ref{og-twist}) while the second is
defined as
$$
{\cal F}^{(s)} = \exp(-2\xi (v_+ \otimes v_+) \phi (\sigma \otimes 1,
1 \otimes \sigma))
$$
where all the terms of the expansion for the symmetric function $\phi$
can be written down explicitely \cite{KUL}. This provides the possibility to
extract the corresponding Lie bialgebra and find that it reveals
the decomposition property (\ref{sum-mu}):
$\left( g_{osp(1|2)},\right.$ $\left. hg_1^{*} + \xi g_2^{*}\right)$.
For the defining relations of ${\cal U}(osp(1|2))$ as
$$
 \left[ h, v_{\pm} \right] = \pm v_{\pm}, \quad \{ v_+ , v_- \}=
 -h/4, \quad X_{\pm} = \pm 4 v_{\pm} v_{\pm}
$$
the corresponding Lie maps (in terms of dual basis) are
\begin{equation}
\label{suplie}
\begin{array}{ll}
\mu _{1}^{*}(\widehat{h},\widehat{X_{\pm}}) 
= -2\widehat{X_{\pm}}, & \quad \mu _{2}^{*}(\widehat{X_{+}},\widehat{h}) 
= 2\widehat{h}, \\
\mu _{1}^{*}(\widehat{h},\widehat{v_{\pm}}) = -\widehat{v_{\pm}},  & \quad
\mu _{2}^{*}(\widehat{X_{+}},\widehat{X_{-}}) = 2\widehat{X_{-}}, \\
\mu _{1}^{*}(\widehat{v_{\pm}},\widehat{v_{\pm}}) = 4\widehat{X_{\pm}}, 
& \quad
\mu _{2}^{*}(\widehat{X_{+}},\widehat{v_{\pm}}) = \widehat{v_{\pm}}, \\
&
\quad \mu _{2}^{*}(\widehat{v_{+}},\widehat{v_{+}}) = 4\widehat{h}, \\
&
\quad \mu _{2}^{*}(\widehat{v_{+}},\widehat{v_{-}}) = 4\widehat{X_{-}},
\end{array}
\end{equation}
Both maps, $\mu _{1}^{*}$ and $\mu _{2}^{*}$, are 2-cocycles of each
other and the second of them is a 2-coboundary $\mu _{2}^{*}
\in B^2(g_1^{*},g_1^{*})$:
$$
\mu _{2}^{*}= \delta \psi,
$$
the 1-cochain $\psi$ on the basic elements looks like
$$
\psi : \left( \widehat{h},\widehat{X_+},\widehat{X_-},
\widehat{v_+},\widehat{v_-}\right) \rightarrow
\left( -\widehat{X_+}, -\widehat{h}, -\widehat{X_-}, \widehat{v_-},
 \widehat{v_-} \right).
$$

The cohomological properties of the involved Lie bialgebras permit to
explain the connection of the Drinfeld-Jimbo (standard) quantization with
the twisting. This is a special case of the general dependence:
Having the first order coboundary deformation of a Lie bialgebra $\left(
g,g_1^{*}\right)$
$$
\left( g,g_1^{*}\right) \longrightarrow \left( g,g_1^{*} + \xi
g_2^{*}\right)
$$
with $\mu_2^{*} \in B^2 \left( \mu_1^{*},\mu_1^{*} \right)$ one can always
get the quantized Lie bialgebra ${\cal A}\left( g,g_2^{*}\right)$ as a limit
of the sequence of quantizations of the
type ${\cal A}\left( g,g_1^{*}\right)$.


\begin{thebibliography}{99}
\bibitem{DRIN}   Drinfeld V. \,G.: ''Quantum groups'', in: {\em Proc. Int.
Congress of Mathematicians, Berkeley, 1986}, {\em 1}, ed. A. V. Gleason
(AMS, Providence, 1987) 798.

\bibitem{JIMB}  Jimbo M.: Lett. Math. Phys. {\em 10} (1985) 63; {\em 11}
 (1986) 247.

\bibitem{ETI}  Etingof P.,  Kazhdan D.: Selecta Math. {\em   2}  (1)
(1996) 1; q-alg/9510020.

\bibitem{D2}  Drinfeld  V. \, G.: Leningrad Math. J. {\em   1}, (1990) 1419,
 DAN USSR {\em   273}, (3) (1983) 531.

\bibitem{GER}  Gerstenhaber M., Giaquinto A.,  Schack S. \ D. in {\it 
Quantum groups. Proceedings in EIMI 1990, Lect. Notes Math.} {\em   1510} 
ed. P. P. Kulish (Springer-Verlag, 1992) 9.

\bibitem{OGIEV}  Ogievetsky O. \, V. in {\em Proc. Winter School Geometry
and Physics, Zidkov, Suppl. Rendiconti cir. Math. Palermo, Serie II -- N 37},
(1993) 185.

\bibitem{KLM} Kulish P. \, P., Lyakhovsky V.\, D., Mudrov A.\, I.: "Extended 
jordanian twists for Lie algebras",  
preprint math.QA/9806014.

\bibitem{ABD} Abdesselam B., Chakrabarti A., Chakrabarti R.: ''General
Construction of Nonstandard $R_h$-matrices as Contraction Limits of $R_q$
-matrices'', preprint q-alg/9706033.

\bibitem{NIJ} Nijenhuis A., Richardson R.: Bull. Amer. Math. Soc. {\em   
72}, (1966) 1.

\bibitem{Lya} Lyakhovsky V.\, D., Mirolyubov A.\, M.: Intern. J. Mod.
Phys. A {\em   12}, (1997) 225.

\bibitem{KUL}  Kulish P. \, P.: "Super-jordanian deformation of the 
orthosymplectic Lie superalgebras", preprint math.QA/9806104
\end{thebibliography}
\end{document}